\documentclass{amsart}

\usepackage[english]{babel} 
\usepackage{amsmath,amsfonts,amssymb,amscd}
\usepackage[pdftex]{graphicx}

\newtheorem{theorem}{Theorem}
\newtheorem{proposition}[theorem]{Proposition}
\newtheorem{corollary}[theorem]{Corollary}

\newtheorem{conjecture}[theorem]{Conjecture}
\newtheorem{definition}[theorem]{Definition}
\newtheorem{example}[theorem]{Example}
\newtheorem{remark}[theorem]{Remark}

\title{A filtration of (q,t)-Catalan numbers}

\author{N. Bergeron}\address[Nantel Bergeron]{Department of Mathematics and Statistics\\ York  University\\
   To\-ron\-to, Ontario M3J 1P3\\ CANADA} \email{bergeron@mathstat.yorku.ca}  \urladdr{http://www.math.yorku.ca/bergeron}

\author{F. Descouens}\address[F. Descouens]{The Fields Institute\\ 222 College Street\\
Toronto, Ontario, M5T 3J1\\ CANADA} \email{francois.descouens@utoronto.ca}
\urladdr{http://www.fields.utoronto.ca/\~{}chohlweg}

\author{M. Zabrocki} \address[Mike Zabrocki]{Department of Mathematics and Statistics\\ York  University\\
    To\-ron\-to, Ontario M3J 1P3\\ CANADA} \email{zabrocki@mathstat.yorku.ca}  \urladdr{http://www.math.yorku.ca/\~{}zabrocki}

 \thanks{This work is supported in part by CRC and NSERC. It is a result of the Algebraic
Combinatorics Seminar at the Fields Institute with the active
participation of Anouk Bergeron-Brlek, Philippe Choquette, and  Muge Taskin. }

\date{}
 
\begin{document}
\maketitle
\def\RR{{\mathbf{R}}}

%%%%%%%%%%%%%%%% TITLE %%%%%%%%%%%%%%%%%%
%\title{A generalization of (q,t)-Catalan and nabla operators} 

%\author{N. Bergeron}
%\address[Nantel Bergeron]
%{Department of Mathematics and Statistics\\ York  University\\ To\-ron\-to, Ontario M3J 1P3\\ CANADA}
%\email{bergeron@mathstat.yorku.ca}
%%\urladdr{http://www.math.yorku.ca/bergeron/}

% \author{F. Descouens}
%\address[F. Descouens]{The Fields Institute\\
%222 College Street\\
%Toronto, Ontario, M5T 3J1\\ CANADA}
% \email{}
%% \urladdr{}

% \author{M. Zabrocki}
% \address[M. Zabrocki]
% {Department of Mathematics and Statistics\\ York  University\\ To\-ron\-to, Ontario M3J 1P3\\ CANADA}
% \email{zabrocki@mathstat.yorku.ca}
%% \urladdr{http://garsia.math.yorku.ca/\~zabrocki/}

%\date{\today}

%\thanks{This work supported in part by CRC and NSERC}

%%%%%%%%%%%%%%%%%% ABSTRACT %%%%%%%%%%%
\begin{abstract}
Using the operator $\nabla$ of F. Bergeron, Garsia, Haiman and Tesler \cite{FGHT}
acting on the $k$-Schur functions \cite{LLM, LM1, LM2} 
indexed by a single column has a coefficient
in the expansion which is an analogue of the $(q,t)$-Catalan number with a level $k$.
When $k$ divides $n$ we conjecture a representation theoretical
model in this case such that the graded dimensions of the module are the coefficients
of the $(q,t)$-Catalan polynomials of level $k$.
When
the parameters  $t$ is set to $1$, the Catalan numbers of level $k$
are shown to count the number of Dyck paths that lie below a certain Dyck path with $q$ counting the area of the path.  
\end{abstract}

%%%%%%%%%%%%%%%%%% INTRO %%%%%%%%%%%%%
\section{Introduction}
In the study of the $(q,t)$-Catalan polynomials, Bergeron, Garsia, Haiman and Tesler \cite{FGHT}
introduced a remarkable operator on symmetric functions, $\nabla$,
to help explain the conjectured graded Frobenius series of the space of the diagonal harmonic
alternants.
The operator $\nabla$ has Macdonald's symmetric functions as eigenfunctions
(see equation \eqref{nablasym} for a definition) and it was a necessary tool for
arriving at a combinatorial formula for the $(q,t)$-Catalan numbers \cite{GarsiaHaglund}.
The original definition of the $(q,t)$-Catalan numbers is equivalent to the
coefficient of the symmetric function $s_{1^n}(X)$ in the expression $\nabla s_{1^n}(X)$.

The authors Lapointe, Lascoux, Morse \cite{LLM} introduced
and Lapointe, Morse \cite{LM1}, \cite{LM2}, \cite{LM3}, \cite{LM4} further developed  
an analogue of the Schur basis of the space of symmetric functions that they called
$k$-Schur functions.  Here the parameter $k \geq 1$ indicates a level of a filtration of the space
of symmetric functions and the parts of the partitions indexing the $k$-Schur functions
are all less than or equal to $k$.  In summary, the 
$k$-Schur functions $\{ s_\lambda^{(k)}(X;t) \}_{\lambda_1 \leq k}$
are the `fundamental' basis of the space linearly spanned by the
elements $\{ s_\lambda (X/(1-t)) \}_{\lambda_1 \leq k}$ where $f(X/(1-t))$
is the symmetric function $f(X)$ with the primitive power sum elements $p_k(X)$
replaced by $p_k(X)/(1-t^k)$.  $k$-Schur functions are a remarkable analogue
of the Schur basis and the Schur functions and the $k$-Schur functions 
are equal when $k \rightarrow \infty$.
In special cases, $k$-Schur functions are equal to Hall-Littlewood symmetric functions,
but in general there is currently no relatively simple definition of these symmetric functions.

Since the $k$-Schur functions are an analogue of the Schur functions, we
decided to consider the action of the operator $\nabla$ on these symmetric functions,
in particular in the case when $k$-Schur functions are indexed by a single column. 
We found that when $\nabla$ acts on the $k$-Schur function $s_{1^n}^{(k)}(X;1/t)$
then this expands positively again in the $k$-Schur functions
$s_\lambda^{(k)}(X;1/t)$.  In fact, our experiments suggest that $\nabla$ acting on the
$k$-Schur functions $s_{\lambda}^{(k)}(X;1/t)$ where $\lambda =(a^b)$ is a rectangle
also expands positively in the $k$-Schur
functions with inverted parameter, 
however it is not true for arbitrary $\lambda$ that $\nabla s_\lambda^{(k)}(X;1/t)$
again lies in the space linearly spanned by the $k$-Schur functions with
inverted parameter (the first failed example is $\nabla s_{2211}^{(4)}(X;1/t)$)

We define a version of the $(q,t)$-Catalan polynomials that includes a level $k$
by setting $C_n^{(k)}(q,t)$ to be the coefficient of $s_{(1^n)}^{(k)}(X;1/t)$ 
in the expression $\nabla s_{1^n}^{(k)}(X;1/t)$ or more simply
(or, more simply, their definition is $ \langle s_{1^n}(X), \nabla s_{1^n}^{(k)}(X;1/t) \rangle$).  Experimental evidence and
special cases lead us to believe that these numbers form a filtration of the
$(q,t)$-Catalan numbers (see Conjectures \ref{conjpositive} and \ref{conjfiltration})
and so we suspect that $C_n^{(k)}(q,t)$ is a $(q,t)$-counting of some subsets of
Dyck paths.  In certain cases, we can provided a combinatorial interpretation of
these polynomials in terms of subsets of Dyck paths.
Because the operator $\nabla$ is an algebra homomorphism at $t=1$, we
are able to give an explicit interpretation of the polynomials $C_n^{(k)}(q,1)$
as the sum over $q$ raised to the area statistic for each Dyck path which lies below
the path which has $k$ steps up followed by $k$ steps over, followed by
$k$ steps up, followed by $k$ steps over, etc., followed finally by $n~mod~k$
steps up and $n~mod~k$ steps over.

It is interesting to remark that $s_{(1^n)}^{(k)}(X;1/t)$ is equal to the modified Hall-Littlewood polynomials\break
$t^{-n(\mu)}\omega Q'_\mu(X;t)$ where $\mu = (k^{n \, \mathrm{div}\, k}, n \, \mathrm{mod} \, k)$  [See Eq.~(\ref{kSchurHL})]. 
From this one can refine our filtration further and define for any partition $\mu\vdash n$ the $\mu$-Catalan number  $C_{\mu}(q,t)$ to be the coefficient of $s_{(1^n)}(X)$ 
in the expression $\nabla t^{-n(\mu)}\omega Q'_\mu(X;t)$. This would would give a filtration of the $(q,t)$-Catalan numbers compatible with the dominance order
of partitions of $n$. All the results and conjecture presented here work in the same way. We do not consider that generality since J.~Haglund and J.~Morse~\cite{HagMorse} 
have comunicated to us an even more refined definition indexed by composition of $n$, 
see Remark~\ref{HM}.

The remainder of this paper is divided into 4 sections.  In
section \ref{sectionbasicdefs}, we discuss necessary definitions.  
In section \ref{sectiongencats},
we introduce the Catalan numbers indexed by a level $k$ and consider special cases
and specializations.  In section \ref{sectionextra}, we briefly consider the analogous filtrations
of the Schr\"oder paths and parking functions.  Finally in the last section we define a
filtration of the space of diagonal harmonic alternants.  We conjecture based on experimental
data that for $k$ dividing $n$ that the graded dimensions of this space  
are given by the polynomials $C_n^{(k)}(q,t)$.

%\end{section}

%%%%%%%%%%%%%%%%%%%%% SECTION 1 %%%%%%%%%%%%%%%%%%%%%%%%%%%%%%%%%%
\section{Basic definitions} \label{Section1} \label{sectionbasicdefs}
%% Subsection on commutative symmetric functions
%%%%%%%%%%%%%%%%%%%%%%%%%%%%%%%%%%%%%%%%%%%%%%%%%
\subsection{Symmetric functions}
 For symmetric functions, we mainly follow the notations of \cite{Macdo}.
Let $X=\{x_1,x_2,\ldots\}$ be a sequence of  variables. The complete homogeneous symmetric function of degree $n$ in the variables $X$ are defined by
  $$
  h_n(X)=\sum_{i_1\le i_2\le\cdots\le i_n} x_{i_1}x_{i_2}\cdots x_{i_n}\,.
  $$
The space of symmetric functions $Sym$ over a field $F$ is the polynomial ring $F[h_1,h_2,\ldots]$, where $h_n=h_n(X)$. This is a graded ring where $\deg(h_n)=n$.
It is convenient to index bases of $Sym$ by partitions which are sequences $\lambda=(\lambda_1,\lambda_2,\ldots,\lambda_k)$ where $\lambda_1\ge\lambda_2\ge\cdots\ge \lambda_k>0$.  The sequence $\lambda$ is a partition of $n$ if $n=\lambda_1+\cdots+\lambda_k $ and its length $\ell(\lambda)$ is $k$. The homogeneous basis can be defined by
$$h_\lambda=h_{\lambda_1}h_{\lambda_2}\cdots h_{\lambda_k}\ .$$
The elementary basis is defined by $e_\lambda=e_{\lambda_1}e_{\lambda_2}\cdots e_{\lambda_k}$, where 
 $e_n$ is defined by the recurrence 
 $$\left \lbrace
 \begin{array}{ccl}
 e_{-k} & =&  h_{-k} = 0  \text{ for } k > 0\ ,\\
 e_0   & = & 1 \ , \\ 
 0  & =  & \sum_{i+j=n} (-1)^i h_ie_j \ .
 \end{array}\right .
$$
 And for any partition $\lambda$ of $n$, the Schur basis can be defined in an algebraic way by
 $$ s_\lambda =\det\left[ h_{\lambda_i+i-j} \right]_{1\le i,j \le n}. $$
 The usual scalar product on the space $Sym$ is defined on the Schur basis by
 \begin{equation}
 \langle s_\lambda\ , \ s_\mu \rangle = \begin{cases} 1 & \text{ if } \lambda = \mu\ , \\ 0 & \text{ otherwise. } \end{cases} 
 \end{equation}
%\end{subsection}
%% Definition of Macdonald polynomials and Nabla
%%%%%%%%%%%%%%%%%%%%%%%%%%%%%%%%%%%%%% 
\subsection{Macdonald polynomials and Hall-Littlewood functions}
For any partition $\lambda$, we denote by $\lambda^\prime$ the conjugate partition of $\lambda$. The usual normalization constant $n(\lambda)$ is defined by
\begin{equation}
n(\lambda) = \sum_{i \ge 1} (i-1)\lambda_i \ .
\end{equation}
Let us now recall some basic definitions on Macdonald polynomials. The modified Macdonald polynomials $\widetilde{H}_\lambda(X;q,t)$ are defined by
\begin{equation}
\widetilde{H}_\lambda(X;q,t) = t^{n(\lambda)} J_\lambda\left (X;q,\frac{1}{t}\right )\ ,
\end{equation} 
where $J_\lambda(X;q,t)$ is the integral version of Macdonald polynomials defined in VI.8 of \cite{Macdo}. The modified Hall-Littlewood polynomials can be obtained as a specialization of the Macdonald polynomials
$$Q'_\lambda(X;t)=\widetilde{H}_\lambda(X;0,t)\ .$$
The linear operator $\nabla$ introduced in \cite{FGHT} is defined by
\begin{equation}\label{nablasym}
\nabla \widetilde{H}_\lambda(X;q,t) = t^{n(\lambda)}q^{n(\lambda^\prime)} 
\widetilde{H}_\lambda(X;q,t) \ . 
\end{equation}
There exist a long list of conjectures about the action of $\nabla$ on different bases of symmetric functions. For many of them, recent work in this area has developed combinatorial models (proved or conjectural \cite{GarsiaHaglund, GarsiaHaglund2, HHLRU, HagCatalan, HagS, HagLoehr, Loehr,LW}) which explains the different properties.
%\end{subsection}
%% k-Schur Functions
%%%%%%%%%%%%%%%%%
\subsection{$k$-Schur functions}
%{\bf 
%
% MIKE... can you write this}
 % k-Schur functions are triangularly related to the Schur functions
 % and satisfy the Pieri rule.  But what about the t?
 The $k$-Schur functions $s_\lambda^{(k)}(X;t)$ of Lapointe, Lascoux, Morse (see \cite{LLM, LM1, LM2, LM3}) are the fundamental basis of the space 
 ${\mathcal L}\{ Q^\prime_\lambda(X;t) \text{ with } \lambda_1 \leq k\}$,
where ${\mathcal L}$ represents the vector space linear span
of the elements.\\
We are interested in this short note only 
in the explicit definition of the elements $s_{1^n}^{(k)}(X;t)$. 
Let $\mu$ be the partition defined by 
\begin{equation}
\mu = (k^{n \, \mathrm{div}\, k}, n \, \mathrm{mod} \, k) \ .
\end{equation}
For these symmetric functions, we simply define them to be
\begin{equation}\label{kSchurHL}
s_{1^n}^{(k)}(X;t) = t^{n(\mu)}\omega \left ( Q^\prime_\mu\left(X;\frac{1}{t}\right)\right )\ .
\end{equation}
This definition permits us to give the explicit expansion of $s_{1^n}^{(k)}(X)$, for the special cases $k > n/2$
\begin{equation}\label{kschurcolumn}
s_{1^n}^{(k)}(X;t) = s_{1^n}(X)\, +\, ts_{21^{n-2}}(X)\, +\,  \ldots \, + \, t^{n-k}s_{2^{n-k}1^{2k-n}}(X) \ .
\end{equation}
For a complete definition of the $k$-Schur functions with the parameter
$t$, we refer the reader to the references \cite{LLM, LM1}.  Note that these two
references provide two different definitions which are conjectured to be equivalent.
In the case of the indexing partition equal to $1^n$ we can show that they are
both equal to \eqref{kschurcolumn}.

\section{Generalizations of $(q,t)$-Catalan numbers} \label{sectiongencats}
%We study a special case of Conjecture \ref{ConjectureNablakSchur}, where $\lambda$ is a column partition. 
The $(q,t)$-Catalan numbers $C_n(q,t)$ defined in \cite{GarsiaHaiman}, are related to the operator $\nabla$ applied to an elementary symmetric function $e_n(X)$. As defined in the previous section, the $k$-Schur functions indexed by column partitions are a generalization of these elementaries functions which are equal to $e_n(X)$, for $k\ge n$. Hence, a natural way to obtain filtrations of $(q,t)$-Catalan numbers is to replace in this picture, the functions $e_n(X)$ by the $k$-Schur functions. With this process, we obtain new polynomials in $q$ and $t$ with positive coefficients, which are smaller than the usual $(q,t)$-Catalan numbers. By specializing $q=1$ and $t=1$ in these filtrations, we obtain different generalizations of Catalan numbers than those given in \cite{HP}.

%% The usual (q,t)-catalan numbers
%%%%%%%%%%%%%%%%%%%%%%%%%%%%%%%%%%
\subsection{$(q,t)$-Catalan numbers}
Let first recall the definition and the combinatorial interpretation for the $(q,t)$-Catalan numbers in terms of Dyck paths \cite{GarsiaHaglund, GarsiaHaglund2, GarsiaHaiman,HHLRU, HagCatalan}.
%% definition
\begin{definition}[(q,t)-Catalan numbers]\label{DefqtCatalan}
The $(q,t)$-Catalan numbers are the polynomials in the parameters $q$ and $t$ defined by 
\begin{equation}
C_n(q,t) = \langle \nabla e_n(X) \ , s_{1^n}(X) \rangle \ ,
\end{equation}
where $\langle \ , \ \rangle$ is the usual scalar product on symmetric functions.
\end{definition}
These polynomials are in $\mathbb{N}\lbrack q,t \rbrack$. Their specialization at $t=1$ and $q=1$ gives the usual Catalan numbers $C_n$
\begin{equation*}
C_n(1,1)=C_n \ .
\end{equation*}
The $(q,t)$-Catalan numbers are symmetric in the variables $q$ and $t$, i.e. $C_n(q,t) = C_n(t,q)$. The maximum degree in these parameters are 
\begin{equation*}
\deg_q (C_n(q,t)) = \deg_t (C_n(q,t)) = \left ( {}_{2}^{n} \right ) \ .
\end{equation*} 
\begin{example}\label{Catalan6} For $n=6$, the $(q,t)$-Catalan number $C_6(q,t)$ can be represented by an array, where the entry $(i,j)$ corresponds to the coefficient of $q^i\ t^{\left ( {}_{2}^{n} \right )-j}$.
\scriptsize
$$q^i\quad 
\begin{array}{|cccccccccccccccc}
  &   &   &   &   &   &   &   &   &   &   &   &   &   &   & 1  \\
  &   &   &   &   &   &   &   &   &   &   &   &   &   & 1 &    \\
  &   &   &   &   &   &   &   &   &   &   &   &   & 1 & 1 &    \\
  &   &   &   &   &   &   &   &   &   &   &   & 1 & 1 & 1 &    \\
  &   &   &   &   &   &   &   &   &   &   & 1 & 1 & 2 & 1 &    \\
  &   &   &   &   &   &   &   &   &   & 1 & 1 & 2 & 2 & 1 &    \\
  &   &   &   &   &   &   &   &   & 1 & 1 & 2 & 3 & 2 &   &    \\
  &   &   &   &   &   &   &   & 1 & 1 & 2 & 3 & 3 & 1 &   &    \\
  &   &   &   &   &   &   & 1 & 1 & 2 & 3 & 4 & 2 & 1 &   &    \\
  &   &   &   &   &   & 1 & 1 & 2 & 3 & 4 & 3 & 2 &   &   &    \\
  &   &   &   &   & 1 & 1 & 2 & 3 & 4 & 3 & 2 &   &   &   &    \\
  &   &   &   & 1 & 1 & 2 & 3 & 4 & 3 & 2 & 1 &   &   &   &    \\ 
  &   &   & 1 & 1 & 2 & 3 & 3 & 2 & 2 &   &   &   &   &   &    \\
  &   & 1 & 1 & 2 & 2 & 2 & 1 & 1 &   &   &   &   &   &    \\
  & 1 & 1 & 1 & 1 & 1 &   &   &   &   &   &   &   &   &    \\ 
1 &   &   &   &   &   &   &   &   &   &   &   &   &   &    \\\hline 
\end{array}
$$
$$
t^{\left ( {}_{2}^{n} \right )-j}$$
\end{example}
In \cite{HHLRU}, the authors proved that these polynomials have positive coefficients by interpreting them as generating polynomials of Dyck paths with two statistics $area$ for the $t$ and $dinv$ for the $q$. We should also mention that the original combinatorial interpretation, given in \cite{GarsiaHaglund2}, uses the statistics of $area$ and $bounce$ with the two parameters interchanged.
%%
%% Basic stuff on Dyck paths
%%%%%%%%%%%%%%%%%%%%%%%%%%%%
\begin{definition}
A Dyck path of length $n$ is a lattice path from the point $(0,0)$ to the point $(n,n)$ consisting of $n$ north steps and $n$ east steps that never go below the line $y=x$. 
\end{definition}
We denote by $DP_n$, the set of all the Dyck paths of length $n$. Dyck paths of length $n$ are in bijection with sequences $(g_0,\ldots, g_{n-1})$ of $n$ nonnegative integers satisfying the two conditions
\begin{equation}
\left \lbrace
\begin{array}{lcc}
 g_0 = 0 \ ,& \\
 g_{i+1} \le g_i+1\ , & \forall i < n-1 \ .
\end{array}
\right .
\end{equation}
The $i$-th entry $g_i$ of the sequence $g$ corresponds to the number of complete lattice squares between the north step of the $i$-th row of the Dyck path and the diagonal $y=x$. Such sequences are called Dyck sequences. We denote by $DS_n$, the set of all the Dyck sequences of length $n$. From now, we use indifferently Dyck sequences or Dyck paths.
%%
%% Example of Dyck path with its Dyck sequence
%%%%%%%%%%%%%%%%%%%%%%%%%
\begin{example} The Dyck sequence $g=(0,0,1,2,0,1,1,2,3,0)$ corresponds to the Dyck path
\begin{center}
\includegraphics[angle=270, width=3cm]{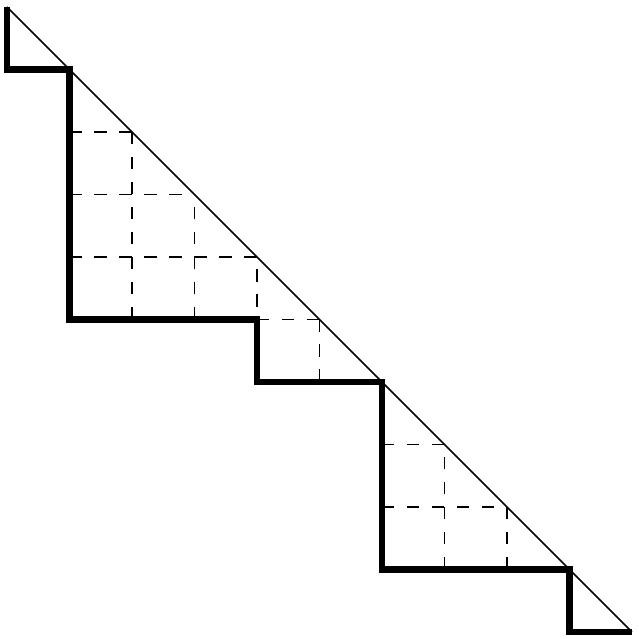}
\end{center}
\end{example}
%%
%% Definition of the two statistics on Dyck paths
%%%%%%%%%%%%%%%%%%%%%%%%%
\begin{definition}
The statistic {\rm area} associated to a Dyck sequence $g$ is defined by
\begin{equation}
\mathrm{area}(g) = \sum_{i=0}^{n-1} g_i \ .
\end{equation}
\end{definition}
On the corresponding Dyck path, this statistic is the number of complete lattice squares between the path and the diagonal $y=x$.
%% Exemple of Area statistic on Dyck path
%% \begin{example} For the Dyck sequence $g$ given in (\ref{ExDyckSequence}), the statistic $area(g)$ is the number of filled squares in the following diagram
%%\begin{center}
%%includegraphics[angle=270, width=4cm]{DyckPathArea.pdf}
%%end{center}
%%\end{example}
\begin{definition}
The statistic {\rm dinv}, which is the number of inversions of a Dyck sequence $g$, is defined by
\begin{equation}
\mathrm{dinv}(g) = \sum_{0\le i < j < n} \chi(g_i-g_j \in \lbrace 0,1\rbrace)\ .	
\end{equation}
\end{definition}
%%
%% Interpretation on the Dyck path
We recall a graphical interpretation of this statistic on Dyck path. Let us call a {\sl north point}, a point where a north step arrives. Two north points give a contribution of 1 in $dinv$, if they are in the same diagonal or if the second point is in the diagonal just below the diagonal of the first one.
%%
%% Graphical exemple of dinv
%%%%%%%%%%%%%%%%%%%%%%%%%%%
\begin{example} The fourth entry of the Dyck sequence $g=(0,0,1,2,0,1,1,2,3,0)$ gives 3 inversions. One inversion comes from the same diagonal and the two others come from the diagonal just below.
\begin{center}
\includegraphics[angle=270, width=3cm]{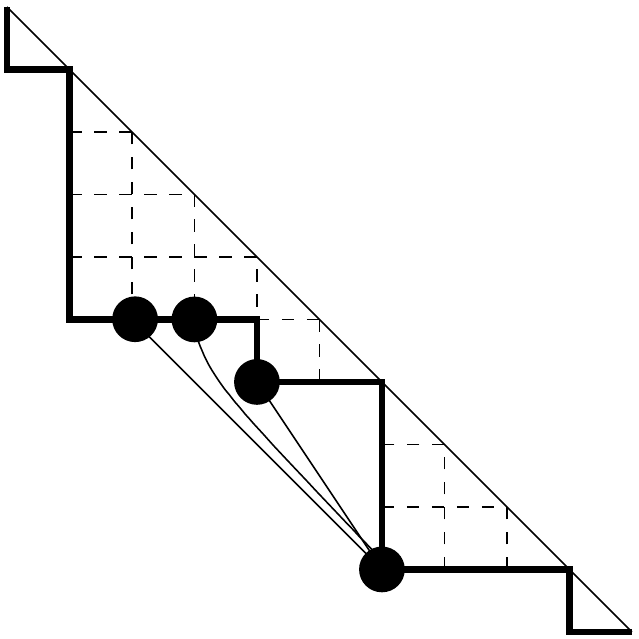}
\end{center}
\end{example}
%% 
%% Combinatorial interpretation 
\begin{theorem}[\cite{GarsiaHaglund2, HHLRU}]
The $(q,t)$-Catalan numbers $C_n(q,t)$ are the generating polynomials of Dyck sequences of size $n$ with the two statistics {\rm area} and {\rm dinv} 
\begin{equation}
C_n(q,t) = \sum_{g \in DS_n} t^{\mathrm{area}(g)}q^{\mathrm{dinv}(g)}\ . 
\end{equation}
\end{theorem}
\begin{example} The $(q,t)$-Catalan number $C_3(q,t)=q^3+q^2t+qt^2+qt+t^3$ can be computing using the following five Dyck paths and the two previous statistics.
\begin{center}
\includegraphics[angle=270, width=7cm]{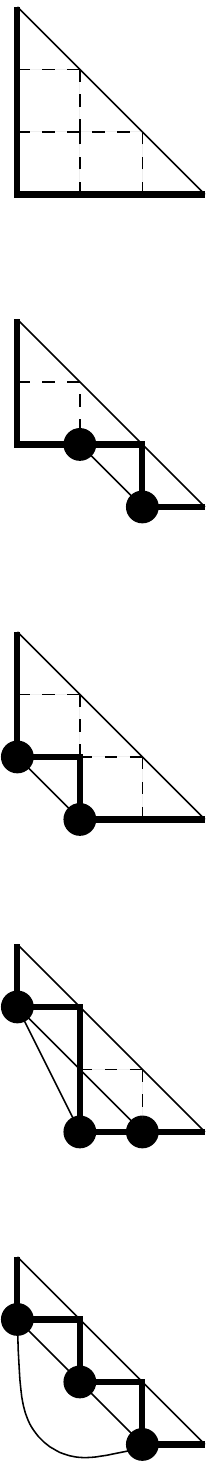}
\end{center}
The black linked points correspond to pairs of points which give a contribution of 1 in the statistic $\mathrm{dinv}$.
\end{example}
%\end{subsection}
%% Definition of the filtration 
%%%%%%%%%%%%%%%%%%%%%%%%%%%%%%%
\subsection{Definition of a filtration of Catalan numbers}
By definition, the $k$-Schur functions indexed by column partitions $(1^n)$ are generalizations of the elementary functions $e_n(X)$ in the space of symmetric functions over $\mathbb{C}(t,q)$. Hence, we can replace elementary functions in Definition \ref{DefqtCatalan} by these $k$-Schur functions, in order to obtain a $k$-level version of $(q,t)$-Catalan numbers.
\begin{definition}
Let $k$ and $n$ be two positive integers. The generalized $(q,t)$-Catalan numbers of level $k$ are defined by
\begin{equation}
C_n^{(k)}(q,t) = \left\langle \ \nabla s^{(k)}_{1^n}\left (X; \frac{1}{t}\right)\ , \ s_{1^n}(X) \right \rangle\ ,
\end{equation}
where $\langle\ , \ \rangle$ is the usual scalar product on symmetric functions.
\end{definition}
\begin{remark}\label{HM}{\rm In fact, $k$-Schur function $s^{(k)}_{1^n}\left (X; \frac{1}{t}\right)$ are special cases of Hall-Littlewood functions (see Equation \ref{kSchurHL}). We can define more general extensions of the $(q,t)$-Catalan numbers indexed by partitions by   
\begin{equation}\label{GenCat}
C_\lambda(q,t) = \left\langle \ \nabla \left (t^{n(\lambda)} \omega\, Q^\prime_\lambda \left (X; \frac{1}{t}\right)\right )\ , \ s_{1^n}(X) \right \rangle \ .
\end{equation}
Recently J. Haglund and J. Morse \cite{HagMorse} have defined generalizations of $(q,t)$-Catalan numbers indexed by compositions using Jing operators on Hall-Littlewood functions and have found the corresponding combinatorial interpretation on Dyck paths. These generalizations are equivalent to our definition (\ref{GenCat}) in the case of partitions. 
}
\end{remark}
%% Positivity conjecture
%%%%%%%%%%%%
\begin{conjecture}[Positivity] \label{conjpositive}
Let $k$ and $n$ be two positive integers. The polynomial $C_n^{(k)}(q,t)$ is in $\mathbb{N}\lbrack q,t \rbrack$.
\end{conjecture}
 In \cite{FGHT}, the authors make the conjecture (originally formulated by Lascoux) that the Schur expansions of the operator $\nabla$ applied on the modified Hall-Littlewood functions $Q^\prime_\lambda(X;t)$ are positive, up to a global sign. The positivity of our generalizations follows directly from this conjecture.
%% Remark on the inversion of the parameter
%%%%%%%%%%%%%%%%%%%%%%%%%%%%%%%%%%%%%%%%%%%
\begin{remark}{\rm
It is important to invert the parameter $t$ inside the $k$-Schur functions in order to obtain polynomials in $q$ and $t$ with integral positive coefficients.}
\end{remark}
%% Increasing properties
%%%%%%%%%%%%%%%%%%%%%%%%
\begin{conjecture}[Filtration]\label{conjfiltration}
Let $n$ be a positive integer. The family of polynomials $\left (C_n^{(k)}(q,t)\right )_{k\ge 1}$ is a filtration of the usual $(q,t)$-Catalan numbers $C_n(q,t)$. More precisely, we have
\begin{equation}
\left \lbrace
\begin{array}{lll}
\forall k\ge 1, & C_n^{(k+1)}(q,t)-C_n^{(k)}(q,t) \in \mathbb{N}[t,q],\\
\forall k\ge n, & C_n^{(k)}(q,t) =C_n(q,t) \ .
\end{array} 
\right .
\end{equation}
\end{conjecture}
{\bf Proof}: 
The first statement is a consequence of Conjecture 3 of \cite{FGHT}. The second statement of the proposition follows immediately from the stability property
\begin{equation}
\forall k \ge n,\ s^{(k)}_{1^n}(X;t) = e_n(X) \ .
\end{equation} 
%
%{\bf We need to use equation \eqref{kschurcolumn} and then use conjecture 3. Nantel will do that later (but soon)}

%% Exemples
%%%%%%%%%%%
\begin{example}\label{Ex6} Using the same conventions as in Example \ref{Catalan6}, the generalized $(q,t)$-Catalan numbers $C_6^{(k)}(q,t)$ are given by the following matrices
\scriptsize
$$
\begin{array}{ccc}
C_6^{(1)} &  C_6^{(2)}& C_6^{(3)} \\\\ 
\begin{array}{|c}
  1 \\\hline 
\end{array} \quad \quad \quad  &
%%
%% level 2
%%
\begin{array}{|cccccccccc}
     &  &  &   &   &   &   &   &   & 1 \\
     &  &  &   &   &   & 1 & 1 & 1 &   \\
     &  &  & 1 & 1 & 1 &   &   &   &   \\
  1  &  &  &   &   &   &   &   &   &   \\\hline 
\end{array} \quad \quad \quad & 
%%
%% level 3
%%
\begin{array}{|cccccccccccccc}
     &  &   &   &   &   &   &   &   &   &   &   & 1 \\
     &  &   &   &   &   &   &   &   &   & 1 & 1 &   \\
     &  &   &   &   &   &   &   & 1 & 1 & 2 & 1 &   \\
     &  &   &   &   &   & 1 & 1 & 2 & 2 &   &   &   \\
     &  &   &   & 1 & 1 & 2 & 1 & 1 &   &   &   &   \\
     &  & 1 & 1 & 1 & 1 &   &   &   &   &   &   &   \\
  1  &  &   &   &   &   &   &   &   &   &   &   &   \\\hline 
\end{array} \quad \quad \quad
\end{array}
$$
$$
%%
%% Level 4
%%
\begin{array}{cc}
\\C_6^{(4)} &  C_6^{(5)} \\
\begin{array}{|ccccccccccccccc}
     &  &   &   &   &   &   &   &   &   &   &   &   & 1 \\
     &  &   &   &   &   &   &   &   &   &   &   & 2 &   \\
     &  &   &   &   &   &   &   &   &   & 1 & 2 &   &   \\
     &  &   &   &   &   &   &   & 1 & 1 & 2 & 1 &   &   \\
     &  &   &   &   &   & 1 & 1 & 2 & 2 &   &   &   &   \\
     &  &   &   & 1 & 1 & 2 & 1 & 1 &   &   &   &   &   \\
     &  & 1 & 1 & 1 & 1 &   &   &   &   &   &   &   &   \\
  1  &  &   &   &   &   &   &   &   &   &   &   &   &   \\\hline 
\end{array} \quad \quad \quad &  
\begin{array}{|cccccccccccccccc}
     &  &   &   &   &   &   &   &   &   &   &   &   &   & 1\\
     &  &   &   &   &   &   &   &   &   &   &   &   & 1 &  \\
     &  &   &   &   &   &   &   &   &   &   &   & 1 & 1 &  \\
     &  &   &   &   &   &   &   &   &   &   & 1 & 1 & 1 &  \\
     &  &   &   &   &   &   &   &   &   & 1 & 2 & 2 &   &  \\
     &  &   &   &   &   &   &   &   & 1 & 2 & 2 &   &   &  \\
     &  &   &   &   &   &   &   & 2 & 2 & 2 & 1 &   &   &  \\
     &  &   &   &   &   & 1 & 2 & 2 & 2 &   &   &   &   &  \\
     &  &   &   & 1 & 1 & 2 & 1 & 1 &   &   &   &   &   &  \\
     &  & 1 & 1 & 1 & 1 &   &   &   &   &   &   &   &   &  \\
  1  &  &   &   &   &   &   &   &   &   &   &   &   &   &  \\\hline 
\end{array}
\end{array}
$$ 
\end{example}
\begin{conjecture}
At the level $k=2$, the generalized $(q,t)$-Catalan numbers statisfy the recursive formula for $n>1$,
$$
\left \lbrace
\begin{array}{ccll}
C^{(2)}_{n+1}(q,t)& = & t^n C^{(2)}_n(q/t,t) &\text{ if $n$ is even}, \\
C^{(2)}_{n+1}(q,t) & = & t^{n} C^{(2)}_n(q,t) + qt^{n-1} C^{(2)}_{n-1}(q,t)&\text{ if $n$ is odd}.\\ 
\end{array}
\right .
$$
\end{conjecture}
%% Example of recursive formula
\begin{example} In the even case, for $n=2$
\begin{equation}
C_3^{(2)}(q,t) = t^2 C_2^{(2)}(q,t) = t^2(t+q/t)= t^3+qt  \ .
\end{equation}
In the odd case, for $n=5$
\begin{eqnarray*}
C_6^{(2)}(q,t) &= &  t^5C_5^{(2)}(q,t) + qt^4C_4^{(2)}(q,t) \\
                         &= & t^5(q^2t^4+qt^7+qt^6+t^{10}) + qt^4(q^2t^2+qt^4+qt^3+t^6) \\
                         &= & q^3t^6+q^2t^9+q^2t^8+q^2t^7+qt^{12}+qt^{11}+qt^{10}+t^{15}
\end{eqnarray*}
\end{example}
\begin{definition}
By specializing $q=1$ and $t=1$ in the generalized $(q,t)$-Catalan numbers $C_n^{(k)}(q,t)$, we define a new filtration $C_n^{(k)}$ of the usual Catalan numbers $C_n$
\begin{equation*}
C_n^{(k)} = C_n^{(k)}(1,1)\ .
\end{equation*}
\end{definition}
%% Exemples a t=q=1
%%%%%%%%%%%%%%%%%%%
\begin{example}
The triangle of the specialization of the generalized $(q,t)$-Catalan at $q=1$ and $t=1$ is
\begin{equation}
\begin{array}{c|cccccc}
n:k & 1 & 2 & 3  & 4  & 5  & 6   \\ \hline
1   & 1 &   &    &    &    &     \\
2   & 1 & 2 &    &    &    &     \\
3   & 1 & 2 & 5  &    &    &     \\
4   & 1 & 4 & 5  & 14 &    &     \\
5   & 1 & 4 & 10 & 14 & 42 &     \\
6   & 1 & 8 & 25 & 28 & 42 & 132     
\end{array}
\end{equation}
\end{example}
The first diagonal below the main diagonal corresponds to the sequence $C_{n-1}$ and the second diagonal below the main one corresponds to the sequence $2\, C_{n-2}$. The others diagonal sequences are unknown in Sloane's integer encyclopedia. But using the combinatorial interpretation at $t=1$, we give an explicit expression for these numbers in the next section.
%% Combinatorial interpretation of qt catalan at t=1
%%%%%%%%%%%%%%%%%%%%%%%%%%%%%%%%%%%%%%%%%%%%%%%%%%%%
\subsection{Combinatorial interpretation at $t=1$}
When the parameter $t$ is specialized at 1 in the generalized $(q,t)$-Catalan numbers $C_n^{(k)}(q,t)$, we are able to give an explicit combinatorial interpretation of these polynomials. This interpretation is based on the fact that $\nabla$ is multiplicative at $t=1$ and was remarked in \cite{FGHT}.  
%% Factorization formula for generalized qt Catalan at t=1
\begin{proposition}
Let $n$ and $k$ be two positive integers. The generalized $(q,t)$-Catalan numbers satisfy the following factorisation formula at $t=1$
\begin{equation}
C_n^{(k)}(q,1)={C_n(q,1)}^{n\ \mathrm{div}\ k}\ C_{(n\ \mathrm{mod}\ k)}(q,1)\ .
\end{equation}
\end{proposition}
{\bf Proof}: By definition, we have 
\begin{eqnarray*}
s_{1^n}^{(k)}(X;1)&=&\omega\left ( H_{(k^{n\ \mathrm{div}\ k},\ n\ \mathrm{mod}\ k)}\left (X; 1 \right ) \right )
                    = \omega\left ( h_{(k^{n\ \mathrm{div}\ k},\ n\ \mathrm{mod}\ k)} (X) \right )\\
                &=& e_{(k^{n\ \mathrm{div}\ k},\ n\ \mathrm{mod}\ k)}(X)\ . 
\end{eqnarray*} 
Now, since the operator $\nabla$ at $t=1$ is multiplicative, we can write
\begin{eqnarray*}
C_n^{(k)}(q,1) &=& \langle \nabla_{t=1}\left ( s_{1^n}^{(k)}(X;1) \right )\ , \ e_n(X) \rangle \\
              &=& \langle \nabla_{t=1}(e_k(X))^{n \ \mathrm{div} \ k}\ 
                           \nabla_{t=1}(e_{n \ \mathrm{mod}\ k}(X))\ , \ e_n(X) \rangle \ .
\end{eqnarray*}
By consideration of degree, the coefficient in $e_n(X)$ in the right part of the scalar product can only be obtained as the product of the coefficient of $e_n(X)$ in $ \nabla_{t=1}(e_k(X))^{n \ \mathrm{div} \ k}$ and in $\nabla_{t=1}(e_{n \ \mathrm{mod}\ k}(X))$ . \hfill $\square$
\begin{remark}
{\rm Using the previous proposition, we obtain an explicit expression for the generalized Catalan numbers $C_n^{(k)}$ in terms of usual Catalan numbers}
\begin{equation}
C_n^{(k)} = \left ( C_{k}\right )^{n \, \mathrm{div}\, k}\, C_{n \, \mathrm{mod} \, k} \ .
\end{equation}
\end{remark}
\begin{corollary}
Let $n$ and $k$ be two positive integers. The combinatorial interpretation of $C_n^{(k)}$ is given by 
\begin{equation}
C_n^{(k)}(q,1) = \sum_{g} q^{\rm{area}(g)}\ , 
\end{equation}
where the sum is taken over the Dyck paths $g$, which are below the Dyck path built with $n\ div\ k$ blocks of $k$ steps up and $k$ steps right and a last block of $n\ mod\ k$ steps up and $n\ mod\ k$ step right. 
\end{corollary}
%% Example of the combinatorial interpretation at t=1
\begin{example}The set of Dyck paths for the combinatorial interpretation of $C_7^{(3)}(q,1)$ are those under the following Dyck path.
\begin{center}
\includegraphics[angle=270, width=2.2cm]{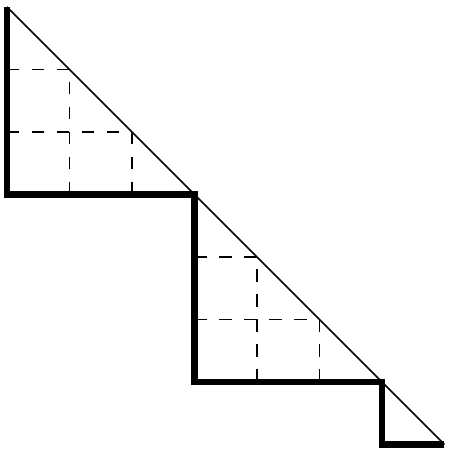}
\end{center}
\end{example}
%
%\end{subsection}
%% Combinatorial interpretation of generalization of (q,t)-Catalan
%%%%%%%%%%%%%%%%%%%%%%%%%%%%%%%%%%%%%%%%%%%%%%%%%%%%%%%%%%%%%%%%%% 
\subsection{Combinatorial interpretation of $C_n^{(k)}(q,t)$}
We find some conjectural combinatorial models for special cases of these generalizations of $(q,t)$-Catalan numbers. We use the combinatorics of configurations of Dyck paths which permits to give some conjectures on a combinatorial model for  $\left \langle \nabla s_\lambda(X)\ , \ s_{1^n}(X)\right \rangle$. It is known that the Schur expansion of $\nabla s_\lambda(X)$ on the Schur basis is always positive, up to a global sign. This sign is interpreted by M. Bousquet-M\'elou using determinants of Catalan numbers. In the special case of $t=1$, a proof using configurations of Dyck paths is given in \cite{Lenart}. More recently, an interpretation using nested Dyck paths is given in \cite{LW} and the following interpretation is mainly based on this work. 

%% Conjectural interpretation of nabla of Schur on Schur basis
%%%%%%%%%%%%%%%%%%%%%%%%%%%%%%%%%%%%%%%%%%%%%%%%%%%%%%%%%%%%%%
\subsubsection{A combinatorial interpretation of $\left \langle \nabla s_\lambda(X)\ , \ s_{1^n}(X)\right \rangle$}
We recall the combinatorial interpretation of the scalar product $\left \langle \nabla s_\lambda(X)\ , \ s_{1^n}(X)\right \rangle$ given in \cite{LW} in terms of nested Dyck paths. For a given partition, we can describe a set of configurations of Dyck paths with two statistics which permits to express the previous scalar product as a generating polynomial. The global sign of these expressions and the characterization of the corresponding configurations of Dyck paths can be computed directly from the partition $\lambda$.  \\\\
%% Removing rim hooks
%%%%%%%%%%%%%%%%%%%%%
For a given partition $\lambda=(\lambda_1, \ldots, \lambda_p)$, we will associate a sequence of $\lambda_1$ nonnegative integers $\tilde{n}(\lambda)=(\tilde{n}_0,\ldots,\tilde{n}_{\lambda_1})$, called the dissection sequence of $\lambda$. Define the maximal rim-hook of a partition $\mu$ as the skew diagram $\mu / (\mu_2-1, \mu_3-1, \ldots, \mu_{\ell(\mu)}-1 )$.  We consider the tiling of the conjugate partition of $\lambda$ obtained by removing successively the maximal rim-hooks. The entry $\tilde{n}_i$ of $\tilde{n}(\lambda)$ is the length of the maximal rim-hook of the tilling starting in the  $(\lambda_1-i)$-th row and is 0 if no rim-hook starts in the $(\lambda_1-i)$-th row.

\begin{example}\label{Dissection} The dissection sequence corresponding to the partition $\lambda=(53222)$ is $\tilde{n}(53222)=(9,0,0,5,0)$ as described by the picture
\begin{center}
\includegraphics[angle=270, width=2.2cm]{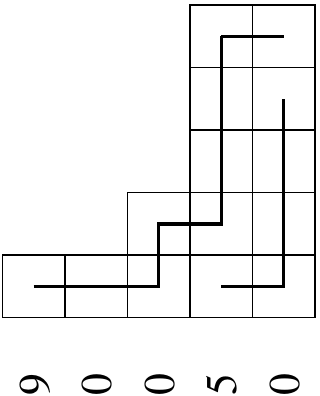}
\end{center}
\end{example}
We define the spin of a partition $\lambda$ by 
\begin{equation}
\mathrm{sp}(\lambda)=\sum_{R} (h(R)-1) \ ,
\end{equation}
where the sum is over all the border rim-hooks of $\lambda^\prime$ and $h(R)$ the height of these ribbons. The sign of a partition $\lambda$, corresponding to the global sign of $\left \langle \nabla s_\lambda(X)\ , \ s_{1^n}(X)\right \rangle $, is defined by 
\begin{equation}
\mathrm{sgn}(\lambda) = (-1)^{\mathrm{sp}(\lambda)} \ .
\end{equation}
We also define the diagonal inversion adjustment by
\begin{equation}
\mathrm{adj}(\lambda) = \sum_{i=0}^{\lambda_1-1}(\lambda_1-1-i)\chi(\tilde{n}_i>0) = \sum_{i=0}^{\lambda_1-1}\lambda_i^\prime\ \chi(\tilde{n}_i>0) \ .
\end{equation}
The adjustment is the sum of the row indices of $\lambda^\prime$ (starting from the top of the diagram of the partition) where a border rim-hook starts. 
%% Example 
%%%%%%%
\begin{example} For the partition $\lambda=(53222)$, the spin and the sign are
\begin{equation}
{\rm sp}(\lambda)=4+1 = 5 \ \text{and consequently}\ {\rm sgn}(\lambda)=-1 \ .
\end{equation} 
In this case, the adjustment is 
\begin{equation}
\mathrm{adj}(\lambda)=1 + 4 = 5 \ .
\end{equation}
\end{example}
%% A nested Dyck path
%%%%%%%%%%%%%%%%%%%%%  
\begin{definition}
Let $\lambda=(\lambda_1,\ldots, \lambda_p)$ be a partition of dissection sequence $\tilde{n}(\lambda)=(\tilde{n}_0, \ldots, \tilde{n}_{\lambda_1-1})$. Let $\Pi=(\pi_0, \ldots, \pi_{\lambda_1-1})$ be a sequence of Dyck paths $\pi_i$ of length $\tilde{n}_i$ from $(i,i)$ to $(i+\tilde{n}_i, i+\tilde{n}_i)$. If $\tilde{n}_i$ is equal to 0, $\pi_i$ is a degenerate Dyck path consisting in a single vertex at $(i,i)$. The sequence $\Pi$ is a nested Dyck path for the partition $\lambda$, if for all $i\not = j$, no edge or vertex of $\pi_i$ coincides with any edge or vertex of $\pi_j$. 
\end{definition}
We denote by $NDP_\lambda$ the set of all the nested Dyck paths for the partition $\lambda$.
%% Example of Nested Dyck path
%%%%%%%%%%%%%%%%%%%%%%%%%%%%%%%
\begin{example}\label{NestedDyckPath} A nested Dyck path of $NDP_{(53222)}$ corresponding to the dissection sequence $\tilde{n}(53222)=(9,0,0,5,0)$.
 \begin{center}
\includegraphics[angle=270, width=2.5cm]{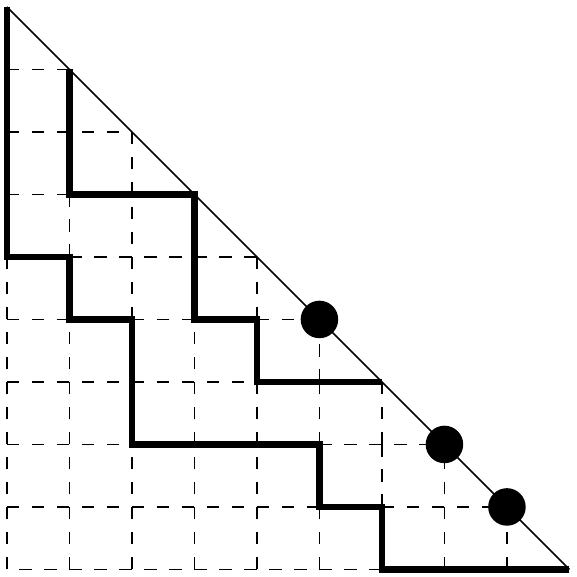}
\end{center}
\end{example}
The encoding of Dyck paths using Dyck sequences can be extended to nested Dyck paths.
Let $\Pi=(\pi_o,\ldots,\pi_{l-1})$ a nested Dyck path. The nested Dyck configuration corresponding to $\Pi$ is an $l$-tuple of words $G=(g^{(0)}, \ldots, g^{(l-1)})$, where $g^{(i)}$ is the Dyck sequence encoding the Dyck path $\pi_i$. The indexing of the letters in these Dyck sequences are chosen to match the alignment of paths in the picture. In the following, we use indifferently nested Dyck paths and nested Dyck configurations.  
%% Example of nested Dyck paths
%%%%%%%%%%%%%%%%%%%%%%%%%%%%%%%
\begin{example} The nested Dyck path of Example \ref{NestedDyckPath} corresponds to the following nested Dyck configuration
\begin{equation}
G = \left ( 
\begin{array}{cccccccccc} 
g^{(0)} : & 0 & 1 & 2 & 2 & 2 & 3 & 4 & 3 & 3  \\
g^{(1)} : & \cdot & \times & \cdot & \cdot & \cdot & \cdot & \cdot & \cdot & \cdot \\
g^{(2)} : & \cdot & \cdot & \times & \cdot & \cdot & \cdot & \cdot & \cdot & \cdot \\
g^{(3)} : & \cdot & \cdot & \cdot &   0   &   1   &   1   &   0   &   1   & \cdot  \\
g^{(4)} : & \cdot & \cdot & \cdot & \cdot & \times & \cdot & \cdot & \cdot & \cdot 
\end{array}
\right ) \ .
\end{equation}
\end{example}
The statistic $area$ and $dinv$ can be extended to nested Dyck paths. The area of a nested Dyck path $G$, written $\overline{area}$, is the sum of the areas of the Dyck paths of $G$:
\begin{equation}
\overline{\mathrm{area}}(G) = \sum_{i=0}^{l-1} {\rm area}(g_i) = \sum_{i=0}^{l-1} \sum_{i\le j < i+n_i}g_j^{(i)} \ .
\end{equation} 
The diagonal inversion statistic for a nested Dyck path $G$, written $\overline{dinv}$, is defined by
\begin{equation}
\begin{array}{ccl}
\overline{\mathrm{dinv}}(G)& = & \mathrm{adj}(\lambda) + \sum_{a,b,u,v} \chi\left( g_a^{(u)}-g_b^{(v)}=1\right)\chi(a\le b)  \\ & & + \sum_{a,b,u,v}
\chi\left( g_a^{(u)}-g_b^{(v)} = 0 \right)\chi( (a<b) \ \text{or} \ (a=b \ \text{and} \ u<v))\ .
\end{array}
\end{equation}
The $\overline{dinv}$ of a nested Dyck path $G=(g^{(0)},\ldots,g^{(l)})$ corresponds to the sum of the $dinv$ of each Dyck path $g^{(i)}$ plus the number of pairs of points coming from different $g^{(i)}$'s which form a inversion. A pair of points which form an inversion and which are in the same row are just counted one time.
 
\begin{example} The statistics $\mathrm{\overline{area}}(G)$ and $\mathrm{\overline{dinv}}(G)$ for the nested Dyck path $G$ of Example \ref{NestedDyckPath} are
\begin{equation}
\mathrm{\overline{area}}(G)=24 \quad \text{and} \quad \mathrm{\overline{dinv}}(G)=37 \ .
\end{equation} 
\end{example}
One of the main conjectures of \cite{LW} gives the following expression for the coefficient of $\nabla s_\lambda(X)$ on the Schur function $s_{1^n}(X)$
\begin{equation}\label{ConjLW}
\langle \nabla s_\lambda(X)  \ , \  s_{1^n}(X)\rangle=\mathrm{sgn}(\lambda) \sum_{G \in NDP_\lambda} q^{\mathrm{\overline{area}}(G)}t^{\mathrm{\overline{dinv}}(G)}\ .
\end{equation}
%% Exemple of the combinatorial interpretation of nabla(s_\lambda)
%%%%%%%%%%%%%%%%%%%%%%%%%%%%%%%%%%%%%%%%%%%%%%%%%%%%%%%%%%%%%%%%%%
\begin{example} For $\lambda=(221)$, we have
\begin{equation}\label{ExNabla}
\langle \nabla s_{221}(X) \ ,\ s_{1^5}(X) \rangle = 
  -\left ( q^6  t^3 + q^5  t^4 + q^4  t^5 + q^3  t^6 \right ) \ .
\end{equation}
The dissection vector of the partition $\lambda=(221)$ is $n=(4,1)$ and 
$\mathrm{adj}(\lambda)= 1$. The combinatorial interpretation of (\ref{ExNabla}) is given by the following four nested Dyck paths where we have linked the pairs of points which give a contribution of 1 in $\overline{\mathrm{dinv}}$.
\begin{center}
\includegraphics[angle=270, width=8cm]{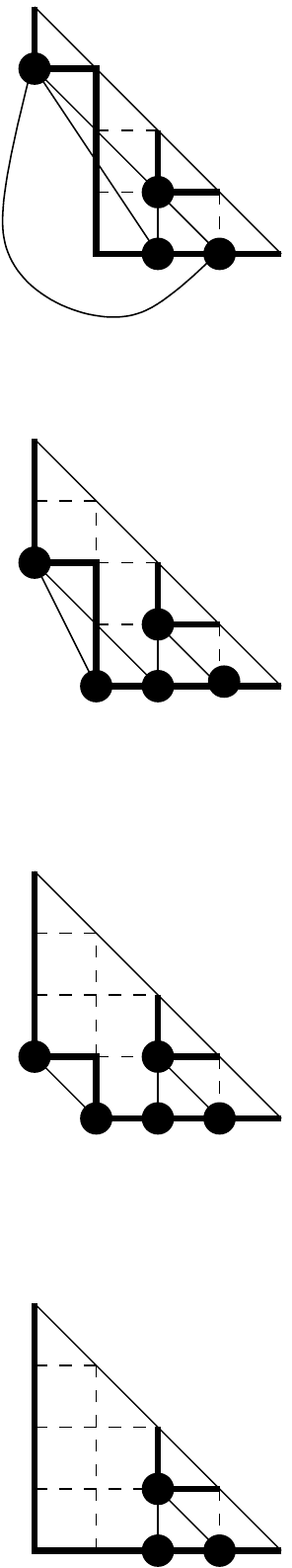}
\end{center}
\end{example}

%\end{subsubsection}
%% Combinatorial interpretation of our filtration
%%%%%%%%%%%%%%%%%%%%%%%%%%%%%%%%%%%%%%%%%%%%%%%%%
\subsubsection{Combinatorial interpretation for the filtration in some special cases}
We give an explicit combinatorial interpretation of the generalizations of $(q,t)$-Catalan numbers in the cases of the level $k=n-1$ and $k=n-2$ using the combinatorial materials given in the previous section. The goal is to find bijections between sets of nested Dyck paths and usual Dyck paths. These bijections have to be compatible with the statistics $(area,\ dinv)$ and $(\overline{area},\ \overline{dinv})$ in order to explain why the terms are canceling in the right way, giving at the end, a polynomial with only positive coefficients.\\\\
%% Combinatorial interpretation for level n-1
%%%%%%%%%%%%%%%%%%%%%%%%%%%%%%%%%%%%%%%%%%%%%
{\bf For level $n-1$}
\\\\
For the level $k=n-1$, we have an explicit characterization of the Schur functions which appear in the $k$-Schur functions we are considering. Using Equation \eqref{kschurcolumn}, we have that 
\begin{equation}\label{kSchurn-1}
s_{1^n}^{(n-1)}(X;t) = s_{1^n}(X) + t\ s_{21^{n-2}}(X) \ .
\end{equation}
%theorem
\begin{conjecture}[Combinatorial interpretation for $k=n-1$]
Let $DP^{(1,1)}_n$ denotes  the set of Dyck paths which go through the lattice point $(1,1)$. The generalized $(q,t)$-Catalan numbers of level $(n-1)$ are given by 
\begin{equation}\label{CatalanLeveln-1}
C_n^{(n-1)}(q,t)=\sum_{g \in DP_n^{(1,1)}} q^{{\rm area}(g)}t^{{\rm dinv}(g)} \ .
\end{equation}
\end{conjecture}
%begin proof of theorem
{\bf Proof (based on Conjecture (\ref{ConjLW}))}: As $s_{1^n}(X)=e_n(X)$, we know that 
\begin{equation}
\langle \nabla s_{1^n}(X)\ , \ s_{1^n}(X)\rangle = \sum_{g \in DP_n} q^{{\rm area}(g)}t^{{\rm dinv}(g)}\ ,
\end{equation}
where the sum is over all the Dyck paths of length $n$.\\
 Let now compute the combinatorial interpretation of $\langle \nabla s_{21^{n-2}}(X)\ , \ s_{1^n}(X)\rangle$ in terms of nested Dyck paths. The dissection vector of the partition $(21^{n-2})$ is $\tilde{n}(21^{n-2})=(n,0)$, as described by the following picture
%% Picture dissection
%%%%%%%%%%%%%%%%%%%%%
\begin{center}
\includegraphics[angle=270, width=3cm]{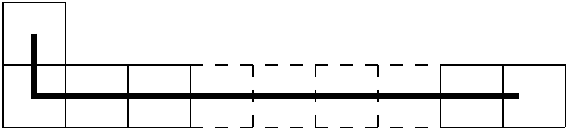}
\end{center}
This implies that the nested Dyck paths corresponding to the partition $(21^{n-2})$ are the sequences of two non intersecting Dyck paths $G=(g^{(0)}, g^{(1)})$, such that
\begin{itemize}
\item[-]Dyck path $g^{(0)}$ is a Dyck path of length $n$ avoiding the lattice point $(1,1)$,
\item[-]Dyck path $g^{(1)}$ is reduced to the degenerated Dyck path of size 0 at the lattice point $(1,1)$. 
\end{itemize}
Hence, we have
\begin{equation}\label{SecondTerm}
\langle \nabla s_{21^{n-2}}(X)\ , \ s_{1^n}(X)\rangle = - \sum_{G \in NDP_{21^{n-2}}} q^{\overline{{\rm area}}(G)}t^{\overline{{\rm dinv}}(G)}\ .
\end{equation}
Let denote by $DP_n^{(1,1)^c}$ the set of all Dyck paths of size $n$ avoiding the lattice point $(1,1)$. Let consider the following bijection $\Phi_n$ defined by
\begin{equation}
\begin{array}{ccccc}
\Phi_n & : & NDP_{21^{n-2}} & \longrightarrow & DP_n^{(1,1)^c} \\
       &   & (g^{(0)}, g^{(1)}) & \longmapsto     & g^{(0)} \ . 
\end{array}
\end{equation}
The compatibility of $\Phi_n$ with the statistics $\mathrm{area}$ and $\mathrm{dinv}$ is given by 
\begin{equation}\label{stat}
\left \lbrace 
\begin{array}{ccl}
\mathrm{dinv}(\Phi_n(g^{(0)}, g^{(1)})) & = & \overline{\mathrm{dinv}}(G) - 1\ , \\
\mathrm{area}(\Phi_n(g^{(0)}, g^{(1)})) & = & \overline{\mathrm{area}}(G) \ .
\end{array} \right .
\end{equation}
%% Proof of the proposition
%%%%%%%%%%%%%%%%%%%%%%%%%%%
In order to prove this compatibility, let $G=(g^{(0)}, g^{(1)})$ be a nested Dyck path of $NDP_{21^{n-2}}$. By definition of $G$, the corresponding Dyck configuration is of the form
\begin{equation}
\left (
\begin{array}{cccccc}
g^{(0)}: & 0 & 1 & g^{(0)}_2 & \cdots & g^{(0)}_{n-1}\\
g^{(1)}: & \cdot & \times & \cdot & \cdots & \cdot
\end{array}
\right ) \ .
\end{equation}
Hence, the Dyck path $g^{(1)}$ always give a contribution of 1 in $\overline{dinv}(G)$.
Using the property of $\Phi_n$ given in (\ref{stat}), Equation (\ref{SecondTerm}) can be rewritten as
\begin{eqnarray}
\langle s_{21^{n-2}}(X)\ , \ s_{1^n}(X)\rangle =& -\ \sum_{G \in NDP_{(21^{n-2})}} q^{\overline{area}(G)} t^{\overline{dinv}(G)} \\
= & -\ \ \sum_{g \in DP_n^{(1,1)^c}}q^{{\rm area}(g)}t^{{\rm dinv}(g)+1} \ . 
\end{eqnarray} 
Hence, we have for generalized $(q,t)$-Catalan of level $(n-1)$
\begin{eqnarray}
C_n^{(n-1)}(q.t)& =& \sum_{g \in DP_n} q^{{\rm area}(g)}t^{{\rm dinv}(g)} - \frac{1}{t}\sum_{g \in DP_n^{(1,1)^c}}q^{{\rm area}(g)}t^{{\rm dinv}(g)+1} \\
& = & \sum_{g \in DP_n^{(1,1)}}q^{{\rm area}(g)}t^{{\rm dinv}(g)} \ .
\end{eqnarray} 
{} \hfill $\square$
\begin{corollary} The generalized Catalan numbers of level $(n-1)$ are given by
\begin{equation}
C_n^{(n-1)}(1,1) = C_{n-1} \ .
\end{equation}
\end{corollary}
The proof is immediate using the combinatorial interpretation given in the previous theorem.\\\\
%% Combinatorial interpretation for level n-2
%%%%%%%%%%%%%%%%%%%%%%%%%%%%%%%%%%%%%%%%%%%%%
{\bf For level $n-2$}
\\\\
In order to give a combinatorial interpretation for generalized $(q,t)$-Catalan numbers of level $(n-2)$, we use the combinatorial interpretation for level $(n-1)$ combined with the combinatorial interpretation of $\langle \nabla s_{221^{n-4}}(X)\ , \ s_{1^n}(X)\rangle$. Using Equation \eqref{kschurcolumn}, we have
\begin{equation}\label{ExprkSchurn-1}
s_{1^n}^{(n-2)}(X;t) = s_{1^n}(X) + t\ s_{21^{n-2}}(X) + t^2\ s_{221^{n-4}}(X) \ .
\end{equation}
%Conjecture level n-2
%%%%%%%%%%%
\begin{conjecture}
Let denote by $DP_n^{(1,1),(3,2)}$ the set of Dyck paths which go through the lattice points $(1,1)$ and $(3,2)$. The generalized $(q,t)$-Catalan numbers of level $(n-2)$ are given by
\begin{equation} 
C_n^{(n-2)}(q,t) = \sum_{g \in DP_n^{(1,1),(3,2)}} t^{\mathrm{dinv(g)}}q^{\mathrm{area}(g)} \ .
\end{equation}
\end{conjecture}
%The proof
{\bf Proof (based on Conjecture (\ref{ConjLW}))}: Let us compute the combinatorial interpretation of $\langle \nabla s_{221^{n-4}}(X)\ , \ s_{1^n}(X)\rangle$. The dissection vector of the partition $(221^{n-4})$ is $\tilde{n}(221^{n-4})=(n-1,1)$ and $adj(221^{n-4})=1$, as described in the following picture
\begin{center}
\includegraphics[angle=270,width=3cm]{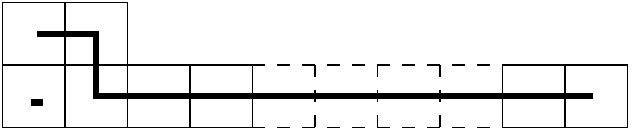}
\end{center}
Hence, a nested Dyck paths $G=(g^{(0)}, g^{(1)})$ is a couple of non intersecting Dyck path $(g^{(0)}, g^{(1)})$ satisfying 
\begin{itemize}
\item[-] $g^{(0)}$ is a Dyck path of length $n-1$ avoiding the lattice point $(2,1)$,
\item[-] $g^{(1)}$ is the unique Dyck path of length 1 starting from the lattice point $(1,1)$.
\end{itemize}
Hence, we have 
\begin{equation}
\langle \nabla  s_{221^{n-4}}(X)\ , \ s_{1^n}(X)\rangle = \sum_{G \in NDP_{221^{n-4}}} q^{\overline{\mathrm{area}}(G)} t^{\overline{\mathrm{dinv}}(G)}\ .
\end{equation}
Let denote by $DP_n^{(1,1),(3,2)^c}$ the set of Dyck paths which go through the lattice point $(1,1)$ and avoid the lattice point $(2,1)$. Let consider the following bijection $\Psi_n$ defined by 
\begin{equation}
\begin{array}{cccc}
\Psi_n : &   NDP_{221^{n-4}}  & \longrightarrow & DP_n^{(1,1),(3,2)^c} \\ 
         & G=(g^{(0)}, g^{(1)}) &   \longmapsto   &    g^{(1)}\cdot g^{(0)} \ , 
\end{array}
\end{equation}
where $g^{(1)}\cdot g^{(0)}$ is the Dyck path of length $n$ obtained by concatenation of $g^{(1)}$ and $g^{(0)}$. 
The compatibility of $\Psi_n$ with the statistics $\mathrm{area}$ and $\mathrm{dinv}$ is given by
\begin{equation}\label{StatisticPsi}
\left \lbrace
\begin{array}{ccl}
\mathrm{dinv}(\Psi_n(g^{(0)}, g^{(1)})) & = &\overline{\mathrm{dinv}}(G)-2\ ,\\
\mathrm{area}(\Psi_n(g^{(0)}, g^{(1)})) & = & \overline{\mathrm{area}}(G) \ .
\end{array}\right .
\end{equation}  
In order to prove this compatibility, let $G=(g^{(0)}, g^{(1)})$ be a nested Dyck path in $NDP_{221^{n-4}}$. The corresponding Dyck configuration is of the form
\begin{equation}
G = \left (
\begin{array}{ccccccc}
g^{(0)} : & 0 & 1 & 2 & g^{(0)}_3 & \cdots & g^{(0)}_{n-2}\\
g^{(1)} : & \cdot & 0  & \cdot & \cdot & \cdots & \cdot
\end{array} \right )
\end{equation} 
The zero of $g^{(1)}$ give a contribution of $2$ in $\overline{dinv}(G)$. By definition of the statistic $dinv$ of a Dyck configuration, we have
\begin{equation}\label{eq1}
\overline{\mathrm{dinv}}(G) = \mathrm{adj}(221^{n-4}) + 2 + \mathrm{dinv}\left(g^{(0)}\right) = 3 + \mathrm{dinv}\left( g^{(0)}\right) \ .
\end{equation}
The concatenation of Dyck paths $g^{(1)}\cdot g^{(0)}$ corresponds to the following Dyck sequence
\begin{equation}
g^{(1)}\cdot g^{(0)} = \left(g^{(1)}_1=0,\ 0,\ 1,\ 2,\ g^{(0)}_3,\ \cdots,g^{(0)}_{n-2}\right) \ .
\end{equation}
The first $0$ gives now a contribution of 1 to $dinv\left( g^{(1)}\cdot g^{(0)}\right)$. Hence, 
\begin{equation}\label{eq2}
\mathrm{dinv}\left( \Psi_n\left(g^{(0)}, g^{(1)}\right)\right) = 
\mathrm{dinv}\left( g^{(1)}\cdot g^{(0)}\right ) = 1 + \mathrm{dinv}\left( g^{(0)}\right)\ .
\end{equation}
Finally, by combining (\ref{eq1}) and (\ref{eq2}), we have $dinv\left( \Psi_n\left(g^{(0)}, g^{(1)}\right)\right)=\overline{dinv}(G) - 2$.\\
Using Equation (\ref{StatisticPsi}), we have
\begin{equation}
\langle \nabla s_{221^{n-4}}(X)\ , \ s_{1^n}(X)\rangle = - \ \sum_{g\in DP_n^{(1,1),(3,2)^c}}q^{\mathrm{area(g)}} t^{\mathrm{dinv(g)} + 2} \ .
\end{equation}
Hence, using Expression (\ref{ExprkSchurn-1}) of $k$-Schur functions when $k=n-2$ and the combinatorial interpretation for level $n-1$, we obtain
\begin{eqnarray*}
\ \ \left \langle \nabla s^{(n-2)}_{1^n}\left(X;\frac{1}{t} \right) \ , \ s_{1^n}(X)\right \rangle  & = & \sum_{g \in DP_n^{(1,1)}} q^{\mathrm{area(g)}}  t^{\mathrm{dinv(g)}}\ - \frac{1}{t^2}
\sum_{g\in DP_n^{(1,1),(3,2)^c}}q^{\mathrm{area(g)}} t^{\mathrm{dinv(g)}+2}  \\
& = & \sum_{g \in DP_n^{(1,1),(3,2)}} t^{\mathrm{dinv(g)}}q^{\mathrm{area}(g)} \ .  
\end{eqnarray*}
{} \hfill $\square$
\begin{corollary}
The generalized Catalan numbers of level $(n-2)$ are given by
\begin{equation}
C_n^{(n-2)}(1,1) = 2\ C_{n-2} \ .
\end{equation}
\end{corollary}
{\bf Proof}: There are two configurations for the first two steps of Dyck paths in $DP_n^{(1,1),(3,2)}$ given in the following picture
\begin{center}
\includegraphics[angle=270,width=5cm]{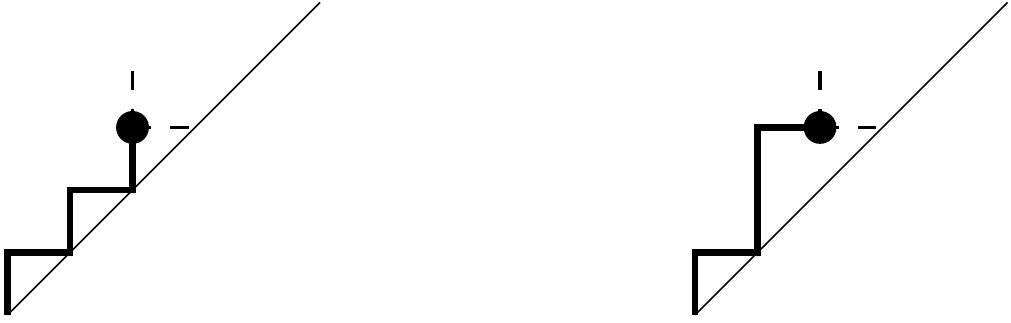}
\end{center}
And it is well known that the number of lattice paths of length $n-3$
starting at the lattice point $(3,2)$ is $C_{n-2}$. Hence the
cardinality of $DP_n^{(1,1),(3,2)}$ is $2\ C_{n-2}$. \hfill $\square$
\\\\
%% Combinatorial interpretation
%%%%%%%%%%%%%%%%%
{\bf For the level $2$}
\\\\
In the special case of $k=2$, we have a conjectural interpretation for the generalized $(q,t)$-Catalan numbers.
\begin{conjecture}
Let $DP_n^{(2)}$ be the set of the Dyck paths which are under the Dyck path given by the sequence  $\left ((10)^{n/2-1} , 1 , (10)^{n/2} , 0\right )$. The generalized $(q,t)$-Catalan number of level k=2 are given by
\begin{equation}
C_n^{(2)}(q,t) = \sum_{g \in DP_n^{(2)}} q^{\mathrm{area}(g)} t^{\mathrm{dinv}(g)} \ .
\end{equation}
\end{conjecture}
%% Algorithm for the other cases
%%%%%%%%%%%%%%%%%%%%%%%%%%%%%%%%
{\bf For other levels $2 < k < n-2$}
\\ \\
For the others levels the problem splits into two different cases. 
For the levels $n/2 < k < n-2$, it exists an algorithm which describe how cancellations are behaving correctly
but the characterizations of the corresponding subsets of Dyck paths are not as nice as for the case of the level $n-1$ and
$n-2$. \\
For the levels $2 < k < n/2$, the coefficients of the $k$-Schur functions indexed by column partitions on the Schur basis are not just monomials in $t$. Hence, the terms are more complicated and are not compatible with the combinatorial interpretation using the same process than before.
%% More ideas on parking functions and Schroeder paths
%%%%%%%%%%%%%%%%%%%%%%%%%%%%%%%%%%%%%%%%%%%%%%%%%%%%%%
\section{Filtration of parking functions and Schr\"{o}der paths} \label{sectionextra}
There exist other interesting polynomials of $\mathbb{N}\lbrack q,t\rbrack$ computed using scalar products involving the operator $\nabla$. For each of these polynomials, a combinatorial model is associated in order to interpret them as generating polynomials with respect with somes statistics. We are mainly interested in the following two examples
\begin{center}
\begin{tabular}{ccc}
$\langle \nabla(e_n(X))\ , \ h_{1^n}(X)\rangle$ & $\longrightarrow$ & $(q,t)$-parking functions, \\
$\langle \nabla(e_n(X))\ , \ e_dh_{n-d}(X) \rangle$ & $\longrightarrow$ & $(q,t)$-Schr\"{o}der paths\ .
\end{tabular}
\end{center}
In order to generalize these combinatorial models, we apply the same idea as in the previous sections, i.e. replacing the elementary functions in $\nabla$ by the $k$-Schur functions indexed by column partitions.\\\\
%
%\subsection{Filtration of the $(q,t)$-parking functions}
The combinatorial model for the $(q,t)$-parking functions have been conjectured by Haglund and Loehr in \cite{HagLoehr,Loehr}. Using $k$-Schur functions, we define a new filtration of these polynomials.
\begin{definition}
Let $k$ and $n$ be two positive integers. The generalized $(q,t)$-parking numbers of level $k$ are defined by 
\begin{equation}
 P_n^{(k)}(q,t) = \left \langle \nabla s_{1^n}^{(k)}\left (X; \frac{1}{t}\right) \ , \ h_{1^n}(X)\right \rangle \ .
\end{equation}
\end{definition}
\begin{example} For $n=3$, the different levels of the filtration are given by
\begin{eqnarray*}
P_3^{(1)}(q,t)  = &  t^3+2t^2+2t+1 \\
P_3^{(2)}(q,t)  = & 2q+2t+2t^2+t^3+qt+1  \\
P_3^{(3)}(q,t)  = & q^3+q^2t+2q^2+qt^2+3qt+2q+t^3+2t^2+2t+1 \ . 
\end{eqnarray*}
\end{example}
The specialization of $t=1$ and $q=1$ in $P_n^{(k)}(q,t)$ gives new sequences of numbers $P_n^{(k)} = (P_k)^{n\ \mathrm{div}\ k} P_{n\ \mathrm{mod}\ k}$, where $P_j=(j+1)^{j-1}$.
$$
\begin{array}{c|cccccc}
n:k &  1  &  2  &  3  &  4  &  5   &  6 \\ \hline
 1  &  1  &     &     &     &    \\
 2  &  2  &   3 &     &     &    \\
 3  &  6  &   9 &  16 &     &     \\
 4  & 24  &  54 &  64 & 125 &     \\
 5  & 120 & 270 & 480 & 625 & 1296 \\
 6  & 720 &2430 &5120 &5625 & 7776 & 16807   
\end{array}
$$
%\begin{proposition}
%The parking functions counted by the numbers $P_n^{(k)}$ are obtained by taking the image of the Dyck paths describing the %generalized $(q,t)$-Catalan numbers of level $k$ by the usual correspondance.
%\end{proposition}
%\begin{example}
%For level $k=n-1$, it corresponds to parking functions with only just one $1$.\\ For level $k=n-2$, it corresponds to parking functions with %just  one $1$ and one $2$ or just one $1$ and two $2$'s.
%\end{example}
The combinatorial model for the $(q,t)$-Schr\"{o}der paths have been conjectured by Egge, Haglund, Kremer and Killpatrick in \cite{EHKK} and proved by Haglund in \cite{HagS}. Using the same kind of idea, we can define filtration of $(q,t)$-Schr\"{o}der paths.
\begin{definition} Let $n,d,k$ be tree positive integers. The generalized $(q,t)$-Schr\"{o}der numbers of level $k$ are defined by
\begin{equation}
\forall d>0, \quad S_{n,d}^{(k)}(q,t) = \left \langle \nabla s_{1^n}^{(k)}\left (X; \frac{1}{t}\right) \ , \ e_d(X)h_{n-d}(X)\right \rangle \ . 
\end{equation} 
\end{definition}
\section{Representation theoretic interpretation of $C_n^{(k)}(q,t)$} \label{sectionreptheory}
The $(q,t)$-Catalan numbers $C_n(q,t)$ are related to the space of
diagonal harmonics $DH_n$ and the $n!$ conjecture on Macdonald
polynomials. Using our generalized $(q,t)$-Catalan numbers
$C_n^{(k)}(q,t)$, we define subspaces $DH_n^{(k)}$ for $k$ dividing $n$ of the space $DH_n$. In the special cases where $k$
divides $n$, we give an explicit algebraic description of these
spaces. We briefly recall some basic statement on the space of
diagonal harmonics and the operator theorem of Haiman \cite{Haiman}. 

%%
%% The Haiman opertors conjecture
%%
\subsection{A generalization of the space of diagonal harmonics $DH_n$}
Let $n$ be a positive integer and $\mathbb{Q}[X_n,Y_n]$ the space of
polynomials over $\mathbb{Q}$ in the two sets of variables $X_n=\lbrace
x_1, x_2, \ldots,x_n \rbrace$ and $Y_n=\lbrace y_1, y_2, \ldots,y_n
\rbrace$. We call bidegree of a polynomial $f(X_n,Y_n)$, the couple of
non-negative integers $(i,j)$ such that $\deg_X(f)=i$ and
$\deg_Y(f)=j$.
\\\\
The symmetric group $S_n$ acts diagonaly on
$\mathbb{Q}[X_n,Y_n]$ by
\begin{equation}
\forall \ f \in \mathbb{Q}[X_n,Y_n],\quad
\sigma.f(x_1,\ldots,x_n,y_1,\ldots,y_n) =
f(x_{\sigma(1)},\ldots,x_{\sigma(n)},y_{\sigma(1)},\ldots,y_{\sigma(n)})\ .
\end{equation}
Let $I$ be the ideal in $\mathbb{Q}[X_n,Y_n]$ generated by all the
$S_n$-invariant polynomials without constant term. Define the
quotient ring
\begin{equation}
  R_n = \mathbb{Q}[X_n,Y_n] / I \ .
\end{equation}
For each $S_n$-invariant polynomials $P(X_n,Y_n)$ of the ideal $I$, the
component of $P$ in $I$ is bihomogeneous in $X_n$ and $Y_n$. Thus, $I$ is
a bihomogeneous ideal. Consequently, $R_n$ has a structure of a doubly
graded ring, i.e.
\begin{equation}
R_n = \bigoplus_{i,j} \left ( R_n \right )_{i,j} \ , 
\end{equation}
where the subspace $\left (R_n\right )_{i,j}$ consists of all images
of homogeneous polynomials of bidegree $(i,j)$. 
\\\\ 
% harmonics point of view
%%%%%%%%%%%%%%%%%%%%%%%%%%%
Let us denote by $\partial x_i$ (resp. $\partial y_i$) the partial
derivative operator with respect to the variable $x_i$
(resp. $y_i$). Define the scalar product
$\langle\ ,\ \rangle_\partial$ on $\mathbb{Q}[X_n,Y_n]$ by
\begin{equation}
\forall f,g \in \mathbb{Q}[X_n,Y_n]\ , \quad \langle f, \ g
\rangle_\partial = f\left (\partial x_1, \ldots, \partial x_n,
\partial y_1,\ldots, y_n \right )g(x_1,\ldots,x_n,y_1,\ldots,\partial y_n)
\ \vert_{X=Y=0} \ .
\end{equation}
For this scalar product the multplication by $x_i$ (resp. $y_i$) is
the adjoint operator of $\partial x_i$ (resp. $\partial y_i$).
%%
%% Definition of the diagonal harmonics 
%%
\begin{definition}
The space $DH_n$ of the diagonal harmonics is defined by 
\begin{equation}
  DH_n=I^{\perp}=\lbrace h \in \mathbb{Q}[X_n,Y_n]\ \vert\ f(\partial x_1,\ldots, \partial x_n,\partial y_1,\ldots, \partial y_n)h = 0\rbrace \ .
\end{equation}
\end{definition}
This definition of the diagonals harmonics is equivalent to the
following caracterization
\begin{equation}
DH_n=\left \lbrace P(X,Y)\in \mathbb{Q}[X_n,Y_n]\ \text{ such that } \sum_{i=1}^n\partial x_i^h\partial y_i^k P \text{ with } h+k>0 \right \rbrace \ .
\end{equation}
The two rings $DH_n$ and $R_n$ are isomorphic and an explicit
isomorphism $\phi : DH_n \longrightarrow R_n$ can be defined by
\begin{equation}
\begin{array}{cccl}
\phi : & DH_n & \longrightarrow & R_n \\
       &  h  & \longmapsto     & \text{the equivalent class of $h$ modulo $I$} \ .
\end{array}
\end{equation}
%%
%% Definition of the alternating diagonals harmonics
%%
In the space $DH_n$, the subspace $DHA_n$ of the alternating harmonics is defined as the diagonals harmonics which are alternating, i.e.
\begin{equation}
DHA_n=\left \lbrace P(X_n,Y_n) \in DH_n \text{ such that } \sigma P(X_n,Y_n) = -P(X_n,Y_n)\ , \forall \sigma \in S_n  \right \rbrace \ .
\end{equation}
\begin{proposition}[\cite{Haiman}]
Let $C_n$ be the $n$-th Catalan number. The dimension of the space $DHA_n$ is given by
\begin{equation}
\dim DHA_n = C_n \ .
\end{equation}
\end{proposition}
The space $DHA_n$ is a bigraded vector space which can be decomposed as 
\begin{equation}
DHA_n = \bigoplus_{i=1}^{\left ( {}_{2}^{n} \right ) } \bigoplus_{j=1}^{\left ( {}_{2}^{n} \right ) } (DHA_n)_{i,j}\ ,
\end{equation}
where $(DHA_n)_{i,j}$ is the space of the polynomials in $DHA_n$ of bidegree $(i,j)$. The Hilbert series of $DHA_n$ is defined by
\begin{equation}
 \mathcal{F}_{DHA_n}(q,t) = \sum_{i=1}^{\infty} \sum_{j=1}^{\infty} t^iq^j \dim (DHA_n)_{i,j} \ .
\end{equation}
\begin{theorem}[\cite{Haiman}]
The $(q,t)$-Catalan numbers are defined by
\begin{equation}
C_n(q,t) = \mathcal{F}_{DHA_n}(q,t) \ .
\end{equation}
\end{theorem}
%%%%%%%%%%%%%%%%%%%%%%%%%%
%% Explanation of the Operator conjecture
%%%%%%%%%%%%%%%%%%%%%%%%%%
\subsection{The operator theorem}
The structure of the space $DH_n$ can be more explicitly, but not entirely, described using the operator theorem given in \cite{Haiman}. The idea is to introduce differential operators $E_k$ which generate the space $DH_n$
only from the Vandermonde determinant of level $n$ in variables $X$ defined by  
%% Definition vandermonde
\begin{equation}
\Delta_n(X_n) = \prod_{1\le i<j \le n}(x_i-x_j) \ . 
\end{equation}
%% Definition operators Ek
These operators $E_k$ are defined for all $k>0$ by
\begin{equation}
E_k = \sum_{i=1}^n y_i \partial x_i^p \ .
\end{equation}
%% The operator conjecture 
\begin{theorem}[\cite{Haiman}]
The space of diagonal harmonics $H_n$ is the smallest space containing $\Delta_n(X_n)$ and closed under the action of the 
operators $E_p$ for all $1 \le p \le n-1$ and the operators $\partial x_i$ for all $1 \le i \le n$. We write this statement using
the following notation
\begin{equation}
DH_n = \mathcal{L}_{E_1,\ldots,E_{n},\partial x_1,\ldots, \partial x_n}(\Delta_n(X_n)) \ .
\end{equation}
\end{theorem}
%% Remarks on the symetry of our generalization of the operators conjecture
If we consider the operators $F_p$ obtained by interchanging $X_n$ and $Y_n$ in $E_p$, the space $DH_n$ can also be described 
as the smallest space containing $\Delta_n(Y_n)$ and closed under the action of the 
operators $F_p$ for all $1 \le p \le n-1$ and the operators $\partial y_i$ for all $1 \le i \le n$. In that sense, 
the operator conjecture is symmetric. Our generalization of the operators conjecture
is not symmetric because our generalization of $(q,t)$-Catalan numbers are not symmetric and the space $DH_n^{(k)}$ 
does not contain the Vandermonde determinant $\Delta_n(Y_n)$ in variables $Y_n$.
%% Special case for the alternant space
\begin{corollary}[\cite{Haiman}]
The space of the alternants $DHA_n$ is the smallest space containing $\Delta_n(X_n)$ and closed under the action of the operators
$E_p$ for all $1 \le p \le n-1$, i.e.
\begin{equation}
DHA_n = \mathcal{L}_{E_1,\ldots,E_{n}}(\Delta_n(X_n)) \ .
\end{equation}
\end{corollary}
\subsection{Special case when $k$ divides $n$}
\begin{conjecture}\label{conjecturekdivn}
Let $k$ and $n$ be two integers such that $k$ divides $n$ and $d=n/k$. Let us define the space $DHA_n^{(k)}$ by
\begin{equation}
DHA_n^{(k)} = \mathcal{L}_{E_d,E_{d+1},\ldots, E_n}(\Delta_n(X_n)) \ .
\end{equation}
The Hilbert series of $DHA_n^{(k)}$ is given by
$
\mathcal{F}_{DHA_n^{(k)}}(q,t) = C_n^{(k)}(q,t) \ .
$
\end{conjecture}
\begin{example} For $n=8$ and $k=4$, we have the triangle corresponding to $C_8^{(4)}(q,t)$ is
\scriptsize
$$q^i\quad 
\begin{array}{|cccccccccccccccccccccccccc}
  && & & & & & & & & & & & & & & & & & & & & & &1\\
  && & & & & & & & & & & & & & & & & & & & &1&1\\
  && & & & & & & & & & & & & & & & & & &1&1&2&1\\
  && & & & & & & & & & & & & & & & &1&1&2&3&2&1\\
  && & & & & & & & & & & & & & &1&1&2&3&5&3&1\\
  && & & & & & & & & & & & &1&1&2&3&5&6&4&2\\
  && & & & & & & & & & &1&1&2&3&5&6&\boxed{7}&3&1\\
  && & & & & & & & &1&1&2&3&5&6&6&5&2&1\\
  && & & & & & &1&1&2&3&5&5&6&4&3&1\\
  && & & & &1&1&2&3&4&4&4&3&1&1\\
  && & &1&1&2&2&3&2&2&1&1\\
  &&1&1&1&1&1&1 \\
1               \\ \hline
\end{array}
$$
$$ t^{\left ( {}_{2}^{n} \right )-j}$$ \normalsize The boxed entry of
coordinates $(19,7)$ corresponds to the subspace of $DHA_8^{(4)}$ of
bidegree $t^{28-19+1}q^{7-1}=t^{10}q^6$ with dimension 7. And
$$ \mathrm{rank}\left \lbrace 
\begin{array}{cccc}
E_{732222}.\Delta_n(X_n), & E_{642222}.\Delta_n(X_n), & E_{633222}.\Delta_n(X_n), & E_{552222}.\Delta_n(X_n), \\
E_{543222}.\Delta_n(X_n), & E_{533322}.\Delta_n(X_n), & E_{444222}.\Delta_n(X_n), & E_{443322}.\Delta_n(X_n), \\
E_{433332}.\Delta_n(X_n), & E_{333333}.\Delta_n(X_n)
\end{array}
\right \rbrace = 7\ .
$$
\end{example}
This conjecture has been verified up to $n=8$.
%%
%% Extension of the previous conjcture at the whole diagonal harmonic space
%%
%The filtration $AH_n^{(k)}$ given in Conjecture \ref{conjecturekdivn}
%can be extented as a filtration $H_n^{(k)}$ of the whole space of
%diagonal harmonics $H_n$.
%\begin{conjecture}
%Let $n$ and $k$ be two integers such that $k$ divides $n$ and $d=n/k$. Let us define the space $H_n^{(k)}$ by
%\begin{equation}
% H_n^{(k)} = \mathcal{L}_{\partial_1,\ldots,\partial_n,E_d,E_{d+1},\ldots, E_n}(\Delta_n(X_n)) \ .
%\end{equation}
% The Hilbert series of $H_n^{(k)}$ is given by
%\begin{equation}
%\mathcal{F}_{H_n^{(k)}}(q,t) = C_n^{(k)}(q,t) \ .
%\end{equation}
%\end{conjecture}

%% The Bibliography
%%%%%%%%%%%%%%%%%%%  

\end{document}